\documentclass[11pt]{amsart}
\usepackage{ amsmath, amsthm, amsfonts, hyperref, graphicx, ifpdf}
\usepackage{mathrsfs,epsfig}
\usepackage[dvipsnames,usenames]{color}
\hypersetup{
   unicode=false,          
   pdftoolbar=true,        
   pdfmenubar=true,        
   pdffitwindow=false,     
   pdfstartview={},    
   pdftitle={},    
   pdfauthor={},     
   pdfsubject={},   
   pdfcreator={},   
   pdfproducer={}, 
   pdfkeywords={}, 
   pdfnewwindow=true,      
   colorlinks=true,       
   linkcolor=blue,          
   citecolor=blue,        
   filecolor=blue,      
   urlcolor=blue          
}

\newtheorem{theorem}{Theorem}[section]
\newtheorem{cor}[theorem]{Corollary}
\newtheorem{lem}[theorem]{Lemma}
\newtheorem{pro}[theorem]{Proposition}
\newtheorem{remark}[theorem]{Remark}
\newtheorem{Def}[theorem]{Definition}
\theoremstyle{definition}

\DeclareMathOperator{\PSL}{\mathsf{PSL}}%
\DeclareMathOperator{\Isom}{\mathsf{Isom}}%
\newcommand{\CC}{\mathcal{C}}
\newcommand{\Bcal}{\mathcal{B}}

\newcommand{\af}{almost Fuchsian}

\newcommand{\caf}{co-area formula}

\newcommand{\chtm}{cusped hyperbolic 3-manifold}
\newcommand{\eum}{Euclidean metric}
\newcommand{\fg}{fundamental group}
\newcommand{\gp}{good position}
\newcommand{\hg}{hyperbolic geometry}
\newcommand{\hym}{hyperbolic metric}
\newcommand{\hys}{hyperbolic space}
\newcommand{\htm}{hyperbolic 3-manifold}
\newcommand{\igp}{in good position}
\newcommand{\jc}{Jordan curve}
\newcommand{\im}{induced metric}
\newcommand{\is}{incompressible surface}
\newcommand{\kg}{Kleinian group}
\newcommand{\la}{least area}

\newcommand{\mc}{mean curvature}
\newcommand{\maxp}{maximum principle}
\newcommand{\mcr}{maximal cusped region}
\newcommand{\ms}{minimal surface}
\newcommand{\nb}{normal bundle}
\newcommand{\nv}{normal vector}
\newcommand{\pc}{principal curvature}
\newcommand{\qf}{quasi-Fuchsian}

\newcommand{\scc}{simple closed curve}

\newcommand{\Tt}{Teichm\"{u}ller theory}
\newcommand{\tm}{3-manifold}
\newcommand{\tg}{totally geodesic}
\newcommand{\uhs}{upper-half space}
\newcommand{\wrt}{with respect to}

\newcommand{\be}{\begin{equation}}
\newcommand{\ene}{\end{equation}}
\newcommand{\br}{\begin{remark}}
\newcommand{\er}{\end{remark}}
\newcommand{\bl}{\begin{lem}}
\newcommand{\el}{\end{lem}}
\newcommand{\bcor}{\begin{cor}}
\newcommand{\ecor}{\end{cor}}
\newcommand{\bpro}{\begin{pro}}
\newcommand{\epro}{\end{pro}}
\newcommand{\ben}{\begin{enumerate}}
\newcommand{\een}{\end{enumerate}}
\newcommand{\bp}{\begin{proof}}
\newcommand{\ep}{\end{proof}}
\newcommand{\bpo}{\begin{pro}}
\newcommand{\epo}{\end{pro}}
\newcommand{\beq}{\begin{equation*}}
\newcommand{\eeq}{\end{equation*}}
\newcommand{\bear}{\begin{eqnarray}}
\newcommand{\eear}{\end{eqnarray}}
\newcommand{\beqar}{\begin{eqnarray*}}
\newcommand{\eeqar}{\end{eqnarray*}}
\newcommand{\bt}{\begin{theorem}}
\newcommand{\et}{\end{theorem}}

\newcommand{\C}{\mathbb{C}}
\newcommand{\R}{\mathbb{R}}

\renewcommand{\H}{\mathbb{H}}
\newcommand{\Z}{\mathbb{Z}}
\newcommand{\Hscr}{\mathscr{H}}
\newcommand{\Nscr}{\mathscr{N}}

\newcommand{\T}{\mathbf{T}}

\newcommand{\kappabar}{\overline{\kappa}}
\newcommand{\gbar}{\bar{g}}

\newcommand{\sbar}{\overline{s}}

\newcommand{\Gammabar}{\overline{\Gamma}}

\newcommand{\Hcal}{\mathcal{H}}

\newcommand{\U}{\mathbb{U}}

\DeclareMathOperator{\dist}{dist}%
\DeclareMathOperator{\grad}{grad}%
\DeclareMathOperator{\Area}{Area}%
\DeclareMathOperator{\Length}{Length}%

\renewcommand{\Im}{\mathop{\mathrm{Im}}}

\numberwithin{equation}{section}

\allowdisplaybreaks

\def\XXint#1#2#3{{\setbox0=\hbox{$#1{#2#3}{\int}$}
    \vcenter{\hbox{$#2#3$}}\kern-.5\wd0}}

\makeatletter
\def\@citestyle{\m@th\upshape\mdseries}
\def\citeform#1{{\bfseries#1}}
\def\@cite#1#2{{%
  \@citestyle[\citeform{#1}\if@tempswa, #2\fi]}}
\@ifundefined{cite }{%
  \expandafter\let\csname cite \endcsname\cite
  \edef\cite{\@nx\protect\@xp\@nx\csname cite \endcsname}%
}{}
\makeatother

\begin{document}

\title{Closed Minimal Surfaces in Cusped Hyperbolic Three-Manifolds}


\author{Zheng Huang}
\address[Z. ~H.]{Department of Mathematics, The City University of New York, Staten Island, NY 10314, USA}
\address{The Graduate Center, The City University of New York, 365 Fifth Ave., New York, NY 10016, USA}
\email{zheng.huang@csi.cuny.edu}

\author{Biao Wang}
\address[B. ~W.]{Department of Mathematics and Computer Science\\
         Queensborough Community College, The City University of New York\\
         222-05 56th Avenue Bayside, NY 11364, USA\\}
\email{biwang@qcc.cuny.edu}

\date{February 1, 2016}



\begin{abstract}
 Motivated by classical theorems on {\ms} theory in compact {\htm}s, we investigate the questions of existence and
 deformations for {\la} {\ms}s in complete noncompact {\htm} of finite volume. We prove any closed immersed incompressible
 surface can be deformed to a closed immersed {\la} surface within its homotopy class in any {\chtm}. Our techniques highlight
 how special structures of these {\chtm}s prevent any {\la} {\ms} going too deep into the cusped region.

\end{abstract}
\maketitle

\section{Introduction}


\subsection{Minimal surfaces in {\htm}s}
Minimal surfaces are fundamental objects in geometry. In {\tm} theory,
the existence and multiplicity of {\ms}s often offer important geometrical
insight into the structure of the ambient {\tm} (see for instance
\cite{Rub05, Mee06}), they also have important applications in {\Tt},
Lorentzian geometry and many other mathematical fields (see for example
\cite{Rub07, KS07}). By Thurston's geometrization theory, the most common
geometry in a {\tm} is hyperbolic (\cite{Thu80}), and this paper is a part
of a larger goal of studying closed incompressible {\ms}s in {\htm}s.

Before we state our main result, we briefly motivate our effort by
making some historic notes on {\ms} theory in three different types of
{\htm}s, namely, compact {\htm}s, {\qf} manifolds,
and cusped {\htm}s (complete, noncompact, and of finite volume).

Let $M^3$ be a complete Riemannian {\tm} (with or without boundary), and let
$\Sigma$ be a closed surface which is immersed or embedded in $M^3$, then
$\Sigma$ is called a \emph{{\ms}} if its {\mc} vanishes identically,
further we call it \emph{least area} if the area of $\Sigma$ {\wrt} the
{\im} from $M^3$ is no greater than that of any other surface which is
homotopic or isotopic to $\Sigma$ in $M^3$.

A closed surface is called \emph{incompressible} in $M^3$ if the induced map
between the fundamental groups is injective, where we don't require that an
incompressible map to be an embedding. Throughout this paper, we always assume
that a closed {\is} is of genus at least two and is oriented.


In the case when $M^3$ is a closed Riemannian $3$-manifold,
Schoen and Yau (\cite{SY79}) and
Sacks and Uhlenbeck (\cite{SU82}) showed that if $S\subset{}M^3$ is a
closed {\is}, then $S$ is homotopic to an immersed {\la} {\ms} $\Sigma$
in $M^3$. The techniques of \cite{SY79, SU82} extend to the case
$M^3$ is a compact (negatively curved) {\tm} with mean convex boundary
(i.e. $\partial{}M^3$ has non-negative {\mc} {\wrt} the inward normal
vector), then there still exists an immersed {\la} {\ms} $\Sigma$ in any
homootopy class of {\is}s (see \cite{MSY82, HS88}). Note that the
existence of immersed closed surfaces in closed {\htm}s (to start
the minimization process) follows from recent remarkable resolution of
the surface subgroup conjecture by Kahn-Markovic (\cite{KM12}).

Recall that a {\qf} manifold is a complete (of infinite volume)
{\htm} diffeomorphic to the product of a closed surface and $\R$. Since
the convex core of any geometrically finite {\qf} manifold is compact
with mean convex boundary, one finds the existence of closed {\is} of
least area in this class of {\htm}s. In \cite{Uhl83}, Uhlenbeck initiated
a systematic study of the moduli theory of {\ms}s in {\htm}s, where she
also studied a subclass of {\qf} manifolds which we call {\emph{{\af}}.
$M^3$ is called {\af} if it admits a closed {\ms} of {\pc}s less than one
in magnitude. Such a {\ms} is unique and embedded in the {\af} manifold
(see also \cite{FHS83}), and therefore one can study the parameterization
of the moduli of {\af} manifolds via data on the {\ms} (see for instance
\cite{GHW10, HW13, San13}). For the uniqueness and multiplicity questions
of {\ms}s in {\qf} manifolds, or in general {\htm}s, one can refer to
\cite{And83, Wan12, HL12, HW15} and references within.

This paper will address the existence question for immersed closed
incompressible {\ms}s in another important
class of {\htm}s: {\chtm}s. $M^3$ is called a \emph{{\chtm}} if it is a complete non-compact {\htm} of finite volume. There are many examples of this type,
frequently the complements of knots and links in the 3-sphere $\mathbb{S}^3$.
Mostow rigidity theorem (\cite{Mos73}) extends to this class of {\htm}s
by Prasad (\cite{Pra73}), however the techniques used in
\cite{SY79, SU82} to find incompressible {\ms}s do not. It is well-known
that any {\chtm} admits infinitely many immersed
closed {\ms}s (\cite{Rub05}), however, they may not be embedded, nor
incompressible. Using min-max theory, very
recently, Collin, Hauswirth, Mazet and Rosenberg in
\cite[Theorem A]{CHMR14} proved the existence of an embedded
(not necessarily incompressible) compact {\ms} in $M^3$. It has been a
challenge to show the existence of immersed
(or embedded) closed incompressible {\ms} in {\htm}s.

For the rest of the paper, we always assume $M^3$ is an oriented {\chtm}.
\subsection{Main result}
In {\tm} theory, it is a question of basic interest to ask if one can deform
an immersed surface in its homotopy class to some area minimizing surface.
Instead of looking for the existence of an oriented, immersed, closed, incompressible minimal surface in a {\chtm} $M^3$, we aim to prove that
one can deform any immersed closed {\is} into a {\la} {\ms} in its
homotopy class. More specifically, we show:

\bt\label{main}
Let $S$ be a closed orientable surface of genus at least two, which is
immersed in a {\chtm} $M^3$. If $S$ is incompressible, then $S$ is
homotopic to an immersed {\la} {\ms} in $M^3$.
\et

We use relatively elementary tools, taking advantage special structure of
the cusps. Given the {\chtm} $M^3$ and an immersed {\is} $S$, we make one
truncation to obtain a compact {\tm} $M^3(\tau_4)$ of negative curvature
and {\tg} boundary. The location where this truncation takes place is
determined by $M^3$ and $S$ (see Remark ~\ref{howdeep}).
We obtain quantitative estimates on how deep this {\la} {\ms} can reach
into the cusped region of $M^3$ (see Remark ~\ref{howdeep} and
Corollary~\ref{sigma'}). The geometric structures both in the {\uhs} $\H^3$
and the {\chtm} $M^3$ play crucial role in our arguments
to keep the {\la} {\ms} in the region not arbitrarily far into the cusp. We observe that any cusped region
is a topologically solid torus with the core curve removed, and an area minimizing closed {\is} can only
have certain ways to intersect the boundary of a cusped region.

Our techniques easily apply to the case when an embedded {\is} is in
presence in $M^3$, namely, we prove the following statement,
which was originally shown by Collin, Hauswirth, Mazet and Rosenberg:

\bcor[{\cite[Theorem B]{CHMR14}}]\label{main2}
Let $S$ be a closed orientable embedded surface in a {\chtm} $M^3$ which
is not a $2$-sphere or a torus. If $S$ is incompressible and non-separating,
then $S$ is isotopic to an embedded  {\la} {\ms}.
\ecor

Their original argument for this result is to cut further and further into the
cusp(s), and apply the results \cite{HS88} each time to obtain
a sequence of {\la} {\ms}s, then show there is at least one such {\ms} in the
hyperbolic region by applying two forms of the {\maxp}.

Note that a general existence theorem for an immersed closed essential surface in any
{\chtm} was established in \cite{CLR97}. It is very special that there exist
some {\chtm}s which do NOT admit any {\em embedded} closed essential surfaces (\cite{Hat82}).


\subsection{Outline of the proof}
We actually prove the embedded case first, namely, Corollary ~\ref{main2}.
There are essentially two parts for it. First we modify
the {\hym} in $\H^3$ to obtain a submanifold of $M^3$ in the quotient
with sufficiently long cusped regions, and the modified metric
around all boundaries so that the submanifold is a compact negatively
curved manifold with {\tg} boundaries. By results of
\cite{MSY82, HS88}, there is a {\la} {\ms} $\Sigma$ ({\wrt} the new metric,
not the {\hym}) in the homotopy class of a closed {\is}
$S$ in this compact submanifold. The heart of the argument is then to
guarantee it does not drift into infinity of $M^3$. We deploy a
co-area formula (see Lemma ~\ref{coarea}) as our main tool for this.
We can then show that $\Sigma$ is actually contained in the
subregion of the submanifold which is still equipped with the {\hym}.
Hence $\Sigma$ is a {\la} {\ms} {\wrt} the {\hym}.
It is oriented as well since the surface $S$ is non-separating. To prove
our main theorem, we lift an immersed essential surface to an
embedded non-separating and {\is} in a finite cover of $M^3$, where we
can apply prior arguments and take advantage of the {\hg}
of a {\chtm} to show the existence of a {\la} {\ms} in the homotopy
class of any closed immersed {\is}.
\subsection{Organization}
The organization of the paper is as follows: in $\S$\ref{sec:prelim},
we cover necessary background material and fix some
notations; in $\S$\ref{sec:barrier}, we modify the {\uhs} model of
$\H^3$ to set up hemispheres as barriers for {\ms}s
in $\H^3$; in $\S$\ref{sec:truncated space}, we move down to the
{\chtm} $M^3$ and its {\mcr}s. Using the modification in previous
section we obtain a truncated Riemannian {\tm} of negative curvature.
Finally in $\S$\ref{sec:proof of main theorem}, we prove our main result.

\subsection{Acknowledgement}

We would like to thank Richard Canary, Joseph Maher and Alan Reid for helpful discussions. We also thank the support from PSC-CUNY
research awards. Z. H. acknowledges supports from U.S. NSF grants DMS 1107452, 1107263, 1107367
``RNMS: Geometric Structures and Representation varieties" (the GEAR Network) and a grant from the Simons
Foundation (\#359635, Zheng Huang). It was a pleasure to discuss some aspects of this project at Intensive Period
on {\Tt} and {\tm} at Centro De Giorgi, Pisa, Italy, and Workshop on Minimal Surfaces and Hyperbolic Geometry at IMPA, Rio, Brazil. We thank 
the referee for careful reading and helpful suggestions.
\section{Preliminary}\label{sec:prelim}
\subsection{{\kg}s and {\chtm}s}
We will work with the {\uhs} model of the {\hys} $\H^3$, i.e.
\beq
   \H^3=\{(x,y,t)\in\R^3\ |\ t>0\}\ ,
\eeq
equipped with metric
\be\label{hym}
   ds^2=\frac{dx^2+dy^2+dt^2}{t^2}\ .
\ene
The {\hys} $\H^{3}$ has a natural compactification: $\overline{\H^{3}}=\H^{3}\cup\widehat\C$,
where $\widehat\C=\C\cup\{\infty\}$ is the Riemann sphere. The orientation preserving isometry group of
the {\uhs} $\H^3$ is given by $\PSL_{2}(\C)$, which consists of linear fractional transformations that
preserve the {\uhs}.

A (torsion free) discrete subgroup $\Gamma$ of $\PSL_{2}(\C)$ is called a \emph{{\kg}}, and the
quotient space $M^3=\H^{3}/\Gamma$ is a complete {\htm} whose {\fg} $\pi_{1}(M^3)$ is isomorphic to
$\Gamma$. Conversely, if $M^3$ is a complete {\htm}, then there exists a holonomy
$\rho:\pi_{1}(M^3)\to\PSL_{2}(\C)$ such that $\Gamma=\rho(\pi_{1}(M^3))$ is a (torsion free) {\kg} and
$M^3=\H^3/\rho(\pi_{1}(M^3))$.

Mostow-Prasad's Rigidity Theorems imply that hyperbolic volume is a topological invariant for {\htm}s
of finite volume, that is to say, these {\htm}s are completely determined by their {\fg}s. J{\o}rgensen and
Thurston (see \cite[Chapter 5--6]{Thu80}) proved that the set of volumes of orientable {\htm}s is
well ordered and of order type $\omega^{\omega}$. Since any non-orientable {\htm} is double-covered
by an orientable {\htm}, then the set of volumes of all {\htm}s is also well ordered.

Many examples of the {\chtm} come from the complements of hyperbolic knots \cite[Corollary 2.5]{Thu82}
on $\mathbb{S}^3$. In general {\chtm}s can be described as follows (see \cite[Theorem 5.11.1]{Thu80}):

\bt
A {\chtm} is the union of a compact submanifold which is bounded by tori and a finite
collection of horoballs modulo $\Z\oplus\Z$ actions.
\et

By the works of Marden, Thurston, Bonahon (\cite{Mar74, Thu80, Bon86}),
any closed incompressible surface of genus at least two in a {\chtm} is
always geometrically finite, i.e. it's either {\qf} or essential with
accidental parabolics (see also the proof of Theorem 5.3 in \cite{Wu2004}).
It is well-known that some {\chtm}s do not contain any embedded closed
{\is}s (\cite{Hat82}). A fundamental fact about
any {\chtm} is the following property which can be found in for instance survey \cite{AFW2015}:

\bt\label{thm:LERF}
The fundamental group of a cusped hyperbolic three-manifold is LERF, i.e.
locally extended residually finite.
\et

As a corollary of Theorem \ref{thm:LERF}, if $S$ is a closed incompressible
surface (with genus $\geq{}2$) immersed in a {\chtm} $M^3$, then $S$ can
be lifted to an \emph{embedded} nonseparating closed {\is}, in a finite cover
of $M^3$ (see \cite{Sco1978,Sco1985,Lon1988,Mat2002}).


\subsection{Maximal cusps and {\mcr}s}\label{sec:maximal cusps}
In this subsection, we briefly describe the \emph{maximal cusps and {\mcr}s}
of the {\chtm} $M^3$, and they will play important roles in our construction.
For more details, one can go to for instance \cite{Ada05, Mar07}.

Suppose that $M^3$ has been decomposed into a compact component (which is called the
\emph{compact core} of $M^3$) and a finite set of cusps (or ends), each homeomorphic to $T^{2}\times[0,\infty)$,
where $T^2$ represents a torus. Each cusp can be realized geometrically as the image of some horoball
$\Hcal$ in $\H^3$ under the covering map from $\H^3$ to $M^3$. If we lift any such cusp to the {\uhs} model
$\H^3$ of the {\hys}, we obtain a parameter family of disjoint horoballs.

Assume first that $M^3$ has exactly one cusp, and we lift it to the corresponding set of disjoint horoballs,
each of which is the image of any other by some group element. Expand the horoballs equivariantly until two
first become tangent. The projection of these expanded horoballs back to $M^3$ is called the
\emph{maximal cusped region} of $M^3$, denoted by $\CC$.

Assume that one such horoball $\Hcal$ is centered about $\infty$. We
may normalize the horoball $\Hcal$ so that $\partial\Hcal$ is a horizontal
plane with Euclidean height \emph{one} above the $xy$-plane.
Thus $\Hcal=\{(x,y,t)\ |\ t\geq{}1\}$. Let $\rho:\pi_{1}(M^3)\to\PSL_{2}(\C)$
be the holonomy of $M^3$. Then $\Gamma=\rho(\pi_{1}(M^3))$ is a
(torsion free) {\kg} with parabolic elements. Let $\Gamma_\infty$ be
the parabolic subgroup of $\Gamma$ which fixes $\infty$, it's then well-known
that $\Gamma_\infty$ is generated by two elements $z\mapsto{}z+\mu$ and
$z\mapsto{}z+\nu$, where $\mu$ and $\nu$ are non-trivial complex numbers
which are not real multiples of each other. Obviously $\Hcal$ is invariant
under $\Gamma_\infty$, and the quotient $\Hcal/\Gamma_\infty$ is just
the {\mcr} $\CC$ of $M^3$ described above. Also
$T^2=\partial\Hcal/\Gamma_\infty$ is a torus.

The fundamental domain of the parabolic group $\Gamma_\infty$ in the
horoball $\Hcal$ is denoted by $A\times[1,\infty)$, where
$A\subset\partial\Hcal$ is a parallelogram spanned by the complex numbers
$\mu$ and $\nu$. It is not hard to see that the Euclidean area of $A$,
which is given by $\Im(\mu\bar\nu)$, is the same
as that of the torus $T^2$.

We may equip the horoball $\Hcal$ with the \emph{warped product metric}
$ds^2=e^{-2\tau}(dx^2+dy^2)+d\tau^2$, by letting $\tau=\log{}t$ for
$t\geq{}1$. Then the metric on the {\mcr} $\CC=T^{2}\times[0,\infty)$
can be written in the form
\be\label{metric-mcr}
   ds^2=e^{-2\tau}\,ds_{\textrm{eucl}}^{2}+d\tau^2\ ,
   \quad \tau\geq{}0\ ,
\ene
where $ds_{\textrm{eucl}}^2$ is the standard flat metric on the torus
$T^2$ induced from that of $\partial\Hcal$.

If $M^3$ has more than one cusp, we define the {\mcr} for each cusp
exactly as above. It's possible that the {\mcr}s in a {\chtm} can
intersect.

Now suppose that the {\chtm} $M^3$ has $k$ cusps, whose {\mcr}s are
denoted by $\CC_{i}=T^{2}_{i}\times[0,\infty)$, $i=1,\dots,k$. Let
$\tau_{0}>0$ be the smallest number such that each {\mcr}
$T^{2}_{i}\times(\tau_{0},\infty)$, $i=1,2,\ldots,k$, is disjoint
from any other {\mcr}s of $M^3$.

For any constant $\tau\geq\tau_{0}$, let $M^{3}(\tau)$ be the compact
subdomain of $M^3$ which is defined as follows:
\be
   M^{3}(\tau)=
   M^{3}-\bigcup_{i=1}^{k}\left(T^{2}_{i}\times(\tau,\infty)\right)\ .
\ene
By this construction, $M^{3}(\tau)$ is a compact submanifold of $M^3$
with concave boundary components {\wrt} the inward normal vectors.

For each $i$ with $1 \le i \le k$, we lift $M^3$ to the {\uhs} model
of the {\hys} $\H^3$ such that one horoball $\Hcal_{i}$ corresponding
to the {\mcr} $\CC_{i}$ is centered at $\infty$ and $\partial\Hcal_{i}$
passes through the point $(0,0,1)$. Suppose that
$\Gamma_{\infty}^{i}$ is the subgroup of $\Gamma$, which is generated by
two elements $z\mapsto{}z+\mu_{i}$ and $z\mapsto{}z+\nu_{i}$, where
$\mu_{i}$ and $\nu_{i}$ are
non-trivial complex numbers that are not real multiples of each other.

Now we may define a constant as follows:
\be\label{maxh}
   L_{0}=\max\big\{e^{\tau_{0}},\,
   |\mu_{1}|+|\nu_{1}|,\,\ldots,\,|\mu_{k}|+|\nu_{k}|\big\}\ > 0.
\ene

\br\label{howdeep}
Note that this constant is independent of $S$. 
If $S$ is embedded in $M^3$, we will prove that the closed incompressible {\la} {\ms} $\Sigma$ in Corollary \ref{main2} is
contained in $M^{3}(\tau_3)$, where $\tau_3=\log(3L_{0})$. If $S$ is only assumed to be immersed in $M^3$, then by Theorem
\ref{thm:LERF}, we may lift $S$ to an embedded {\is} in a finite cover
$N^3$ of $M^3$, which is also a {\chtm}. In this case,
we will show that immersed minimal surface in Theorem \ref{main} is
contained in $M^{3}(\widetilde{\tau}_3)$ for
$\widetilde{\tau}_3=\log(3\widetilde{L}_{0})$, where $\widetilde{L}_{0}$
is defined similarly according to the information of the
cusped regions of $N^3$.
\er
\section{Constructing Barriers in Hyperbolic Three-space}\label{sec:barrier}
In this section we work entirely in the {\hys} $\H^3$ instead of the
quotient {\chtm} $M^3$. Our goal will be to construct hemispheres in
$\H^3$ which can be used as barriers for {\ms}s. To do this, we will first
modify the standard {\hym} on $\H^3$ to get a new metric which
is non-positively curved. This procedure gives us the flexibility we need
to obtain barriers.

\subsection{Modifying the {\hys}}\label{subsec:modified hyperbolic space}

For fixed constants $L_{2}>L_{1}>0$, we define a smooth cut-off function
$\varphi:(0,\infty)\to[0,\infty)$ as follows (see Figure \ref{fig:phi}):
\ben
  \item $\varphi(t)=\frac1t$, if $0<t\leq{}L_{1}$;
  \item $\varphi(t)$ is is strictly decreasing on $[L_{1},L_{2})$, with
        $\varphi(L_{1})=\frac{1}{L_1}$ and $\varphi(L_{2})=0$;
  \item $\varphi(t)\equiv{}0$ if $t\geq{}L_{2}$;
  \item We also require $\varphi$ to satisfy the following inequality:
\be\label{varphi}
   0\leq{}\varphi(t)\leq{}\frac1t\ ,
   \quad\text{for all}\ t>0\ .
\ene
\een

\begin{figure}[htbp]
\begin{center}
  \begin{minipage}[t]{0.45\textwidth}
  \centering
  \includegraphics[scale=0.7]{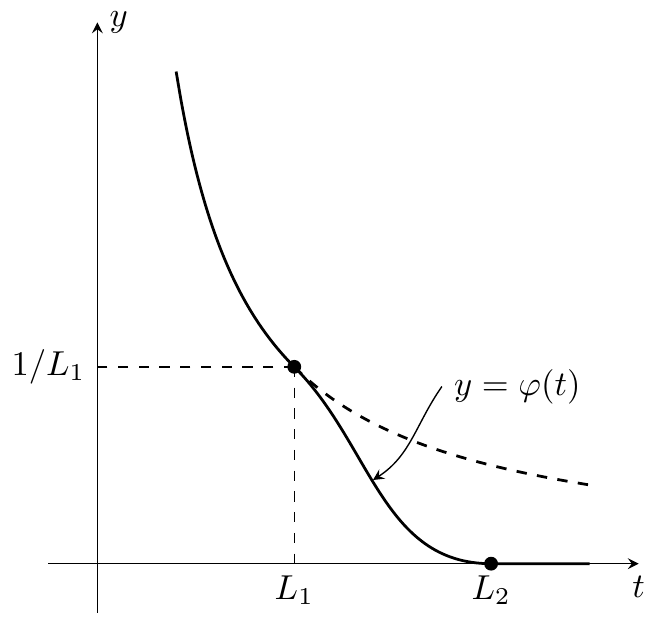}
  \caption{A graph of $\varphi(t)$}\label{fig:phi}
  \end{minipage}%
  \begin{minipage}[t]{0.5\textwidth}
  \centering
  \includegraphics[scale=0.75]{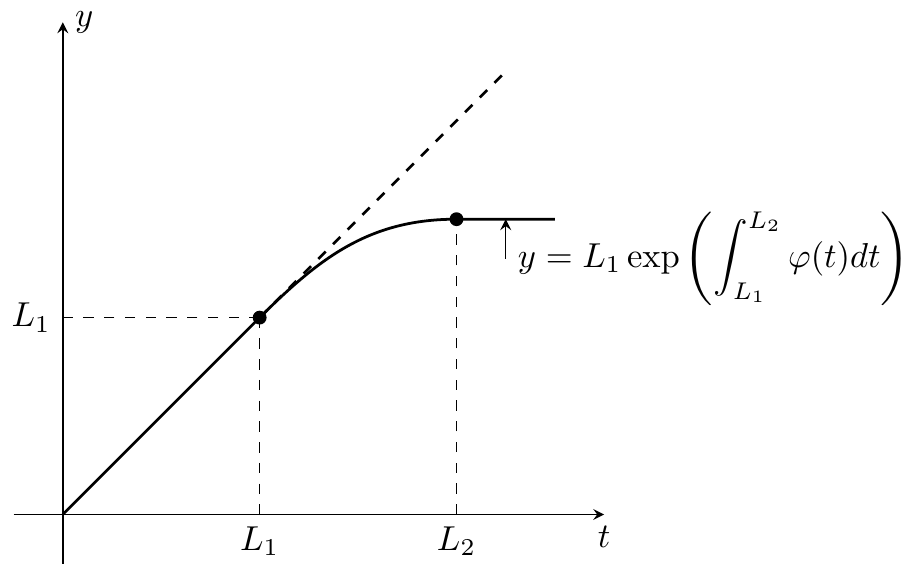}\label{fig:f}
  \caption{A graph of $f(t)$}
  \end{minipage}
\end{center}
\end{figure}

We now define another smooth function $f(t):(0,\infty)\to(0,\infty)$
by solving the following equation:
\be\label{f}
   \frac{f'(t)}{f(t)}=\varphi(t)\ ,
   \quad\text{for all}\ t>0\ .
\ene
And we may require $f(t)$ to satisfy the following (see Figure \ref{fig:phi}) properties:
\ben
\item
$f(t)=t$, if $0<t\leq{}L_{1}$;
\item
$f(t)$ is strictly increasing on the interval $(L_{1},L_{2})$;
\item
$f(t)$ is a constant, if $t\geq{}L_{2}$.
\een

Now we consider an {\uhs} model of the \emph{modified} {\hys} $(\U^{3},\gbar)$,
constructed as follows:
\begin{enumerate}
  \item $\U^3=\R_{+}^{3}=\{(x,y,t)\in\R^3\ |\ t>0\}$,
  \item the new metric is given by
  \be\label{gbar}
  \gbar(x,y,t)=\frac{dx^2+dy^2+dt^2}{(f(t))^2}.
  \ene
\end{enumerate}

Comparing with the standard {\hym} \eqref{hym} on $\H^3$, one sees that $\gbar$ is just
the {\hym} for $t \in (0,L_1]$, and flat beyond $t = L_2$. In fact, we have the following
result, which was not explicitly listed but can be derived from the proof of
\cite[Theorem 4.1]{Zho99}. We include a proof here for the sake of completeness.

\bpro\cite{Zho99}\label{npc}
The {\uhs} $(\U^{3},\gbar)$ is non-positively curved.
\epro

\bp
Recalling from \eqref{gbar}, we may choose a local coordinate system
such that $\gbar_{ij} = \frac{\delta_{ij}}{f(t)^2}$, for $\{i,j\} = \{1,2,3\}$.
We can then workout the Christoffel symbols $\{\bar{\Gamma}_{ij}^k\}$ {\wrt}
this metric $\gbar$ according to the formula:
\beq
\Gammabar_{ij}^k = \frac12 \gbar^{km} (\gbar_{mi,j}+ \gbar_{mj,i} - \gbar_{ij,m}).
\eeq
We find these Christoffel symbols are:

\ben
   \item $\Gammabar^{1}_{13}=\Gammabar^{1}_{31}=
          \Gammabar^{2}_{23}=\Gammabar^{2}_{32}=
          \Gammabar^{3}_{33}=-\frac{f'(t)}{f(t)}$,
   \item $\Gammabar^{3}_{11}=\Gammabar^{3}_{22}=\frac{f'(t)}{f(t)}$, and
   \item all others are equal to $0$.
\een
One can then verify the sectional curvatures of the space $(\U^3,\gbar)$ at a point $(x,y,t)$ are given by
\be
   K_{12}=-(f^{\prime}(t))^2\ ,
   \quad\text{and}\quad
   K_{13}=K_{23}=f''(t)f(t)-(f^{\prime}(t))^2\ .
\ene
Note that, by \eqref{f}, we have
\beq
   \frac{f''(t)f(t)-(f'(t))^2}{f^2(t)}=
   \left(\frac{f'(t)}{f(t)}\right)'=\varphi'(t)\leq{}0\ ,
\quad\text{for all}\ t>0\ .
\eeq
Therefore the space $(\U^3,\gbar)$ is non-positively curved.
\ep

In order to show a convexity statement in Theorem ~\ref{tg}, we need to calculate the {\pc}s of some surfaces immersed in $(\U^3,\gbar)$, if these
surfaces are special {\wrt} to a metric that is conformal to $\gbar$ in $\U^3$. The tool
can be found in the following more general lemma:

\bl[{\cite{Lop13}}]\label{lopez}
For $m \ge 3$, let $(\mathcal{M},g)$ be an m-dimensional Riemannian manifold and let
$\sigma:\mathcal{M}\to\R^{+}$ be a smooth positive function on $\mathcal{M}$. Define
the metric $\gbar=\sigma^{2}{}g$. Let $\iota:S\to{}\mathcal{M}$ be an immersion of an
orientable hypersurface. If $\kappa$ is a {\pc} of $(S,\iota^{*}g)$ {\wrt} the unit normal
vector field $N$, and then
\be\label{kappa}
   \bar\kappa = \frac{\kappa}{\sigma} - \frac{1}{\sigma^2}d\sigma(N)
\ene
is a {\pc} of $(S,\iota^{*}\gbar)$ {\wrt} the unit normal vector field $\overline{N}=N/\sigma$,
and $d\sigma(N)$ is the differential of $\sigma$ along $N$.
\el

By Proposition ~\ref{npc}, we know that the space $(\U^3,\gbar)$ is non-positively curved. We now want to
understand the structure of some special figures in $(\U^3,\gbar)$. This will become important in Theorems
~\ref{tau4} and ~\ref{kg}: we need to construct a submanifold in $M^3$ is of negative curvature and it is a
quotient from a subregion in $\H^3$ by the same {\kg}.

\bt\label{tg}
The subspace $\{(x,y,t)\in\U^3\ |\ 0<t<L_{2}\}$ is a negatively curved space {\rm(}{\wrt}
the metric $\gbar${\rm)}, with a {\tg} boundary $\{(x,y,t)\in\U^3\ |\ t=L_{2}\}$. Furthermore,
any horizontal plane in $(\U^3,\gbar)$ is either convex {\wrt} the upward normal vector
$N=(0,0,1)$, or {\tg}.
\et

\bp
In order to apply Lemma ~\ref{lopez}, on the space $\U^3$, the metric $g$ will
be designated as the Euclidean metric, and the conformal factor
$\sigma(x,y,t) = \frac{1}{f(t)}$, where $f(t)$ is defined previously, and
$\gbar = \frac{g}{f^2(t)}$ is the modified metric on $\U^3$ which is nonpositively
curved by Proposition ~\ref{npc}, and negatively curved in
the subspace $\{(x,y,t)\in\U^3\ |\ 0<t<L_{2}\}$ of $\U^3$.

For any horizontal plane that passes through $(0,0,t)$, its unit {\nv} at the
point $(x,y,t)$ {\wrt} the {\eum} $g$ is
given by $N=\frac{\partial}{\partial{}t}$.

Since
\beq
d\sigma(N)=\grad(1/f(t))\cdot{}N=-\frac{f'(t)}{f^{2}(t)},
\eeq
where $\grad$ is the gradient {\wrt} the Euclidean metric $g$ and $\cdot$ denotes the Euclidean inner
product of vectors, then by \eqref{kappa}, we find the {\pc}s of the plane {\wrt} the new metric $\gbar$
\beq
   \kappabar_{i}(x,y,t)= 0 - f^2(t) (-\frac{f'(t)}{f^{2}(t)}) = f'(t)\ ,\quad i = 1, 2\ .
\eeq
By the construction of the function $f(t)$, we have
\begin{itemize}
  \item $f'(t)>0$ if $0<t<L_{2}$, and
  \item $f'(t)\equiv{}0$ if $t\geq{}L_{2}$.
\end{itemize}
Therefore any horizontal plane through the $(0,0,t)$ is either convex {\wrt} the {\nv} $N=(0,0,1)$
if $0<t<L_{2}$, or {\tg} if $t\geq{}L_{2}$.
\ep

\br\label{vert}
Similarly one can show that any vertical plane is {\tg}, and any vertical straight line is a geodesic {\wrt}
the new metric $\gbar$.
\er
\subsection{Barriers}
The following result guarantees that hemispheres in $(\U^3,\gbar)$
can be used as the barrier surfaces to prevent the {\la} {\ms} $\Sigma$
from entering into each cusped region of $M^3$ too far.

\bt\label{barrier}
For any positive constant $r$, let
\beq
   S_{+}^{2}(r)=\{(x,y,t)\ |\ x^2+y^2+t^2=r^2,\ t>0\}
\eeq
be a hemisphere in $(\U^{3},\gbar)$ with radius $r$. Then
$S_{+}^{2}(r)$ is non-concave {\wrt} the inward {\nv} field, i.e. the {\pc}s
of $S_{+}^{2}(r)$ are nonnegative {\wrt} the inward {\nv} field.
\et

\bp
Let $g$ again denote the standard {\eum} on $\R_{+}^{3}$. At a point $p=\left(x,y,\sqrt{r^2-x^2-y^2}\right)$
on $S_{+}^{2}(r)$, the inward {\nv} field on the hemisphere
$S_{+}^{2}(r)$ {\wrt} the {\eum} $g$ is given by
\beq
   N(p)=\left(-\frac{x}{r},-\frac{y}{r},
           -\frac{\sqrt{r^2-x^2-y^2}}{r}\right)\ .
\eeq
The {\pc}s $\kappa_{1}$ and $\kappa_{2}$ of $S_{+}^{2}(r)\subset(\R_{+}^{3},g)$ {\wrt} the {\nv} $N$ are
identically equal to $\frac1r$.

As in the proof of Theorem ~\ref{tg}, we set $\sigma(x,y,t)=\frac{1}{f(t)}$, where the positive function $f(t)$ is
defined by solving the equation \eqref{f}. Let $\kappabar_{i}$ ($i=1,2$) be the {\pc}s of
$S_{+}^{2}(r)\subset(\U^{3},\gbar)$ at $p$ {\wrt} an orientation $\overline{N}(p)=f\left(\sqrt{r^2-x^2-y^2}\right)N(p)$.

Now we apply \eqref{kappa}, the {\pc}s $\kappabar_{i}$ ($i=1,2$) at $p$ are then given by:
\begin{align*}
   \kappabar_{i}(p)
           &=f\left(\sqrt{r^2-x^2-y^2}\right)\cdot\frac{1}{r}-
             f'\left(\sqrt{r^2-x^2-y^2}\right)\cdot\frac{\sqrt{r^2-x^2-y^2}}{r}\\
           &=\frac{f\left(\sqrt{r^2-x^2-y^2}\right)}{r}
             \left\{1-\varphi\left(\sqrt{r^2-x^2-y^2}\right)
             \sqrt{r^2-x^2-y^2}\right\}\\
           &\geq{}0\ ,
\end{align*}
where we use the property \eqref{varphi}. This completes the proof.
\end{proof}

\section{Truncating Cusped Hyperbolic three-manifold}\label{sec:truncated space}
We want to construct a submanifold in a {\chtm} $M^3$ whose boundary
components are concave {\wrt} the inward normal vectors. The idea is to
remove some horoballs of certain sizes from $\H^3$ in
$\S$\ref{subsec:truncated hyperbolic space}, then modify
the {\hym} in the remaining regions according to previous section, and we
have to of course verify, in $\S$\ref{subsec:Kleinian group}, that the {\kg}
$\Gamma$ of $M^3$ preserves the new metric (otherwise we get a different
{\htm} in the quotient).

\subsection{Truncated {\hys}}\label{subsec:truncated hyperbolic space}
As before we assume that the {\chtm} $M^3$ has $k$ cusps, whose
{\mcr}s are denoted by
$\CC_{i}=T^{2}_{i}\times[0,\infty)$, $i=1,\dots,k$. We also denote $\rho:\pi_{1}(M^3)\to\PSL_{2}(\C)$ as the
holonomy so that $\Gamma=\rho(\pi_{1}(M^3))$ is a {\kg}.

For the $i$-th cusped region $T^{2}_{i}\times[\tau,\infty)$, let
$\Hcal_{i}(\tau)$ be the corresponding horoball
centered at $\infty$, whose boundary is a horizontal plane passing
through the point $(0,0,e^\tau)$, i.e.
\be\label{horo}
\Hcal_{i}(\tau)=\{(x,y,t)\in\H^3\ |\ t\geq{}e^\tau\}.
\ene
In particular, $\Hcal_{i}(0)$ is the corresponding (maximal)
horoball $\Hcal_{i}$ centered at $\infty$. We also
denote $\Hcal_{i}^{\circ}(\tau)$ as the interior of \eqref{horo}.

Recall that $\tau_{0}>0$ is the smallest number such that each {\mcr} $T^{2}_{i}\times(\tau_{0},\infty)$,
$i=1,2,\ldots,k$, is disjoint from any other {\mcr}s of $M^3$. When $\tau\geq\tau_{0}$, the subset
$\Omega(\tau)$ of $\H^3$ is obtained by removing a disjoint collection of open horoballs, namely,
\be\label{tau}
   \Omega(\tau)=\H^3-\bigcup_{i=1}^{k}
   \bigcup_{\gamma\in\Gamma}\gamma \left(\Hcal_{i}^{\circ}(\tau)\right)
\ene
is called a \emph{truncated hyperbolic $3$-space} (see \cite[p.362]{BH99}).

It is clear that $\Omega(\tau)$ is invariant under $\Gamma$, so
\be\label{quotient}
   \Omega(\tau)/\Gamma=M^{3}(\tau)\ .
\ene
We define four constants
\be\label{4c}
   \tau_{j}=\log(j\cdot{}L_{0})\ ,\quad\text{for}\ j=1,2,3,4\ ,
\ene
where the constant $L_{0}$ is defined by \eqref{maxh}. Note that
by this definition \eqref{4c} and by
 \eqref{maxh}, we have $\tau_4 > \tau_3 > \tau_2> \tau_1 \ge \tau_0 > 0$.

We are particularly interested in the subregion $\Omega(\tau_{4})$,
and we define a new metric on it as follows:
\ben
\item
We equip the subregion $\Omega(\tau_{3})$ with the standard {\hym}.
\item
The subregion $\Omega(\tau_{4})\backslash \Omega^{\circ}(\tau_{3})$
(where $\Omega^{\circ}(\tau_{3})$ is the interior of $\Omega(\tau_{3})$)
consists of countably infinitely many disjoint subregions which can
be divided into $k$ families $\Hscr_{1},\ldots,\Hscr_{k}$, such that
each family $\Hscr_{i}$ is the lift of the
cusped subregion $T^{2}_{i}\times[\tau_{3},\tau_4]$.
\een
For an element $U_{i}\in\Hscr_{i}$, we may assume that it can be
described as
\be\label{U}
   U_{i}=\{(x,y,t)\in\H^3\ |\ 3L_{0}\leq{}t\leq{}4L_{0}\}\ .
\ene
We equip the region $U_{i}$ with the new metric
\be\label{f2}
   d\sbar{}^{2}=\frac{dx^2+dy^2+dt^2}{(f(t))^2}\ ,
\ene
where the function $f$ is defined on $[3L_{0},4L_{0}]$ just as in
$\S$\ref{sec:prelim} (i.e. $L_{1}=3L_{0}$ and $L_{2}=4L_{0}$).
Similarly we may define the same new metric on the other elements
in $\Hscr_{i}$, and so on the elements from the other families.

We denote $\gbar$ the new metric on the space $\Omega(\tau_{4})$.
Now we apply Theorem ~\ref{tg} to arrive at the following:

\bt\label{tau4}
The compact space $(\Omega(\tau_{4}),\gbar)$ is a negatively
curved space with (countably infinitely many)
{\tg} boundary components.
\et
\subsection{The {\kg}}\label{subsec:Kleinian group}
The {\kg} $\Gamma$ preserves the hyperbolic metric, but we need to show it also preserves the new metric $\gbar$ on $\Omega(\tau_{4})$. More precisely,

\bt\label{kg}
The group $\Gamma$ is a subgroup of $\Isom(\Omega(\tau_{4}),\gbar)$,
the isometry group of $\Omega(\tau_{4})$ {\wrt} the negatively curved
metric $\gbar$.
\et

\bp
In order not to introduce a different cut-off process, we proceed here with a straightforward (but lengthy) argument.

Let $p$ and $q$ be two points in $\Omega(\tau_{4})$, and we need to show
that $d(p,q)=d(\gamma(p),\gamma(q))$ for any element $\gamma\in\Gamma$,
where $d(\cdot,\cdot)$ denotes the distance function {\wrt} the new
metric $\gbar$. More precisely, let $c$ be the (unique) geodesic
from $p$ to $q$, we shall prove that $\gamma\circ{}c$ is the (unique)
geodesic from $\gamma(p)$ to $\gamma(q)$ for any $\gamma\in\Gamma$.
Moreover we shall prove that $c$ and $\gamma\circ{}c$ have the same
length with respect to the
metric $\gbar$, so $\gamma$ is an isometry of $\Omega(\tau_{4})$ with
respect to the metric $\gbar$ for any $\gamma\in\Gamma$. Therefore $\Gamma$
is a subgroup of $\Isom(\Omega(\tau_{4}),\gbar)$.

By Theorem ~\ref{tau4}, the manifold $(\Omega(\tau_{4}), \gbar)$ is
negatively curved. Then there
is a unique geodesic $c:[0,L]\to(\Omega(\tau_{4}),\gbar)$ parameterized
by arc length, such that
$c(0)=p$ and $c(L)=q$. If the geodesic $c([0,L])$ is totally contained
in $\Omega(\tau_{3})$, we are
done by the definition of the function $f(t)$ (note that $f(t) = t$ for
$t \in (0, 3L_0)$). If $c([0,L])$ is entirely
contained in any component of $\Omega(\tau_{4})-\Omega^{\circ}(\tau_{3})$,
then $f(t)$ is a strictly increasing
function and $\gamma$ preserves the distance.

In general, similar to Corollary 11.34 in \cite{BH99} on page 364,
the geodesic $c$ is expressed as a chain of non-trivial paths
$c_1,\ldots,c_n$, each parameterized by arc length, such that
\ben
  \item each of the paths $c_{i}$ is either a hyperbolic geodesic or else
        its image is contained in one component of
        $\Omega(\tau_{4})-\Omega^{\circ}(\tau_{3})$;
  \item if $c_i$ is a hyperbolic geodesic then the image of $c_{i+1}$ is
        contained in one component of
        $\Omega(\tau_{4})-\Omega^{\circ}(\tau_{3})$, and vice versa.
\een
Suppose that each geodesic segment $c_{i}$ is parameterized by
$c_{i}(s)=c(s)$ for $s\in[s_{i-1},s_{i}]$, where
$0=s_{0}<s_{1}<\dots<s_{n}=L$ is a partition of the interval $[0,L]$.
Then we write $c=c_{1}*c_{2}*\cdots*c_{n}$ in
the sense that $c(s)=c_{i}(s)$ if $s\in[s_{i-1},s_{i}]$. By the above
argument, we have that each curve
$\gamma\circ{}c_{i}:[s_{i-1},s_{i}]\to(\Omega(\tau_{4}),\gbar)$ is a
geodesic for $i=1,\ldots,n$.

We need to show that the curve $\gamma\circ{}c=(\gamma\circ{}c_{1})*\cdots{}*(\gamma\circ{}c_{n})$ is
a geodesic from $\gamma(p)$ to $\gamma(q)$. We will proceed by induction.
To start, $(\gamma\circ{}c_{1})$ is a geodesic segment. Now suppose that $(\gamma\circ{}c_{1})*\cdots{}*(\gamma\circ{}c_{j-1})$ is a geodesic
segment, and $(\gamma\circ{}c_{1})*\cdots{}*(\gamma\circ{}c_{j-1})*(\gamma\circ{}c_{j})$
is not a geodesic segment, then there exists a (unique) geodesic
$c':[0,s_{j}]\to(\Omega(\tau_{4}),\gbar)$ such that $c'(0)=\gamma(p)$ and
$c'(s_{j})=\gamma(c(s_{j}))$, and furthermore the $\gbar$-length of $c'([0,s_{j}])<s_j$. However, $\Gamma$ is a subgroup of $\PSL(2,\C)$, whose
elements are conformal, therefore they preserve the angle. Now three geodesic
segments $(\gamma\circ{}c_{1})*\cdots{}*(\gamma\circ{}c_{j-1})([s_{0},s_{j-1}])$, $\gamma\circ{}c_{j}([s_{j-1},s_{j}])$
and $c'([0,s_{j}])$ would form a geodesic triangle whose sum of its inner
angles is $\geq\pi$. This is a contradiction.

Therefore $\gamma\circ{}c=(\gamma\circ{}c_{1})*\cdots{}*(\gamma\circ{}c_{n})$
is a geodesic segment from $\gamma(p)$ to $\gamma(q)$, and then
$d(\gamma(p),\gamma(q))=L=d(p,q)$.
\ep

As a corollary, we consider the resulting quotient manifold:

\bcor\label{m3tau4}
The manifold $M^{3}(\tau_{4})=\Omega(\tau_{4})/\Gamma$ can be equipped
with a new metric induced from the covering space, still denoted by
$\gbar$, such that $(M^{3}(\tau_{4}),\gbar)$ is a
compact negatively curved {\tm} with {\tg} boundary components.
\ecor

We now make special remarks here on $M^{3}(\tau_{4})$ and its
submanifolds $M^{3}(\tau_{3})$ before we move to the proofs.

\br\label{tau3}
According to the construction of the submanifold $M^{3}(\tau_{4})$,
it is homeomorphic to $M^3$, so its fundamental group
$\pi_{1}(M^{3}(\tau_{4}))$ is also LERF.

By the definition of $f(t)$ in $\S$\ref{subsec:modified hyperbolic space}
and the definition of four constants \eqref{4c}, the modified
metric $\gbar$ restricted to $M^3(\tau_3)$ is the {\hym}. The submanifold
$(M^{3}(\tau_{3}),\gbar)$ is a compact {\htm} whose boundary components are
concave {\wrt} the inward {\nv}s.
\er

\section{Proof of Main results}\label{sec:proof of main theorem}

In $\S$\ref{sec:truncated space} we constructed a submanifold
$M^{3}(\tau_{4})=\Omega(\tau_{4})/\Gamma$ in any {\chtm} $M^3 = \H^3/\Gamma$
with a modified metric $\gbar$ such that $(M^{3}(\tau_{4}),\gbar)$ is a
compact negatively curved {\tm} with mean convex boundary components with
respect to the inward normal vectors.

In this section we may assume that $S$ is an \emph{embedded} closed
incompressible surface with genus $\geq{}2$ contained in
$M^{3}(\tau_{1})\subset{}M^{3}$ until
we begin to prove Theorem~\ref{main} on page \pageref{page:main theorem}.
Now by the argument in \cite{MSY82, HS88}, there exists an
\emph{embedded} closed incompressible {\la} {\ms} $\Sigma$
isotopic to $S$ in $M^{3}(\tau_{4})$ with respect to the modified metric
$\gbar$. In this case, we will say that $\Sigma$ is an embedded
least area minimal surface isotopic to $S$ in $(M^{3}(\tau_{4}),\gbar)$.

\br
If $S$ is only guaranteed to be immersed,
fortunately we then use an additional fact that $\pi_2(M^3) = 0$, and apply
\cite{SU82, HS88} to find the existence of an immersed {\la} surface $\Sigma$
homotopic to $S$ in $(M^{3}(\tau_{4}),\gbar)$.
\er

The key will be showing that the embedded {\ms} $\Sigma$ is
contained in $(M^{3}(\tau_{3}),\gbar)$, a hyperbolic subregion of
$(M^{3}(\tau_{4}),\gbar)$. A key ingredient of the
rest of the argument is that a cusped region has very simple geometry,
and an embedded closed {\is} of the least area can only intersect the region
(if at all) in a predictable way (see Proposition ~\ref{claim1}).

\subsection{Minimal surface intersecting toric region}
As before, we assume that the oriented {\chtm} $M^3$ has $k$ cusps, such
that each maximal cusped region is parametrized by
$\CC_{i}=T^{2}_{i}\times[0,\infty)$ for $i=1,\ldots,k$. Suppose that
$\rho:\pi_{1}(M^3)\to\PSL_{2}(\C)$ is the holonomy so that $\Gamma=\rho(\pi_{1}(M^3))$.

Let $S_{g,n}$ be a surface of genus $g$ with $n$ boundary components
(i.e. a closed genus $g$ surface with $n$ disjoint open disk-type
subdomains removed). If $S_{g,n}$ has negative Euler characteristic, i.e.
$\chi(S_{g,n})<0$, then $S_{g,n}$ must satisfy one of the following
conditions:
\begin{itemize}
  \item If $g\geq{}2$, then $n\geq{}0$.
  \item If $g=1$, then $n\geq{}1$.
  \item If $g=0$, then $n\geq{}3$.
\end{itemize}
It's easy to verify that $\pi_{1}(S_{g,n})$ is non-abelian
in the above three cases.

For a {\scc} $\alpha\subset{}S_{g,n}$, it is said to be \emph{essential}
if no component of $S_{g,n}\setminus\alpha$ is a disk, and it is said
to be \emph{non-peripheral} if no component of $S_{g,n}\setminus\alpha$
is an annulus. We denote $\T'$ as a solid torus with a core curve removed
and we state a simple lemma to be used later:

\bl\label{euler}
Let $S_{g,n}$ be a surface of negative Euler characteristic embedded in
$\T'$ such that $\partial{}S_{g,n}\subset\partial{}\T'$ if $n\geq{}1$.
Then there exists at least one essential simple closed
curve $\alpha\subset{}S_{g,n}$ such that $\alpha$ bounds a
disk $D\subset{}\T'$.
\el

\bp
We prove that the homomorphism $\pi_{1}(S_{g,n})\to\pi_{1}(\T')$ induced by
the embedding $S_{g,n}\to\T'$ can't be injective.
In fact, if any essential simple closed curve $\alpha$ in $S_{g,n}$ is
non-contractible in $\T'$, then the embedding induces an injection
between {\fg}s $\pi_{1}(S_{g,n})\to\pi_{1}(\T')=\Z\oplus\Z$.
Thus $\pi_{1}(S_{g,n})$ is isomorphic to a subgroup of $\Z\oplus\Z$.
But this is impossible, since $\pi_{1}(S_{g,n})$ is non-abelian whereas
$\Z\oplus\Z$ is abelian. Therefore the kernel of the homomorphism $\pi_{1}(S_{g,n})\to\pi_{1}(\T')$ is nontrivial.

Next we shall prove the existence of an essential simple closed curve
in $S_{g,n}$, which represents an element in the kernel
of the homomorphism $\pi_{1}(S_{g,n})\to\pi_{1}(\T')$. Since $\chi(S_{g,n})<0$,
it contains a separating simple closed curve $\alpha$, in other words homologically trivial, so it is a commutator, and 
must be contained in the kernel of the homomorphism $\pi_{1}(S_{g,n})\to\pi_{1}(\T')$.
\ep

Recall that there are four constants only depending on $M^3$:
$\tau_j = \log(jL_0)$, for $j = 1,2,3,4$. Suppose that the embedded
closed {\is} $S$ is contained in $M^{3}(\tau_{1})$.

We consider a compact submanifold $M^{3}(\tau_{4})$ of $M^3$, equipped
with the new metric $\gbar$ (see Corollary ~\ref{m3tau4}), so that
$(M^{3}(\tau_{4}),\gbar)$ is a compact negatively curved {\tm} whose
boundary components are all {\tg}.
By above arguments, we have an embedded closed incompressible {\la} {\ms}
$\Sigma$ isotopic to $S$ in $(M^{3}(\tau_{4}),\gbar)$. Furthermore,
$\Sigma$ is disjoint from $T_{i}^{2}\times\{\tau_{4}\}$ for $i=1,\ldots,k$,
since each boundary component of $(M^{3}(\tau_{4}),\gbar)$ is {\tg}.

We need to show the {\la} {\ms} $\Sigma$ is contained in
$(M^{3}(\tau_3),\gbar)$, that is to say, $\Sigma$ is a minimal surface
with respect to the hyperbolic metric.
If $\Sigma$ does not intersect with $T^{2}_i\times\{\tau_2\}$, then we
are done (since it can not be contained
entirely in the cusped region). Therefore we can just assume that $\Sigma\cap\left(T^{2}_{i}\times[0,\tau_{2}]\right)$ is non-empty.
We are interested in how $\Sigma$ intersects with the region $T^{2}_{i}\times[0,\tau_{4})$.

\bpro\label{claim1}
Each component of $\Sigma\cap(T^{2}_{i}\times[0,\tau_{4}))$ is either
a minimal disk whose boundary is a null-homotopic {\jc} in
$T^{2}_{i}\times\{0\}$, or a minimal annulus whose boundary consists
of essential {\jc}s in $T^{2}_{i}\times\{0\}$.

Moreover each component of $\Sigma\cap(T^{2}_{i}\times[0,\tau_{4}))$
is boundary compressible, i.e. each component can be isotoped into
$T^{2}_{i}\times\{0\}$ such that the isotopy fixes the boundary of
the component.
\epro

\bp
Let $\Sigma'$ be a component of
$\Sigma\cap(T^{2}_{i}\times[0,\tau_{4}))$. Since
$(M^{3}(\tau_{4}),\gbar)$ is a compact negatively curved {\tm} with
{\tg} boundary components, so the {\la} {\ms}
$\Sigma$ is disjoint from its boundary.
Therefore the boundary of $\Sigma'$ is contained in $T^{2}_{i}\times\{0\}$.
Since $\Sigma$ is incompressible while $T^2_i$
is a torus, we have very few cases to consider:
\vskip 0.1in
\textbf{Case I}: Suppose that $\Sigma'$ is a surface of negative
Euler characteristic. We suppose that $\Sigma'$ is homeomorphic to
the surface $S_{g,n}$ with $\chi(S_{g,n})=2-2g-n<0$.

Firstly, no non-peripheral essential curves in $\Sigma'$ are null
homotopic in the region $T_{i}^{2}\times[0,\tau_{4})$.
Otherwise such a curve is also a non-peripheral essential curve in
$\Sigma$. This is impossible since the surface $\Sigma$ is incompressible.

Secondly no peripheral essential curves in $\Sigma'$ are null homotopic
in the region $T_{i}^{2}\times[0,\tau_{4}]$ either. Otherwise,
let $\alpha$ be a boundary component of $\Sigma'$ which is null homotopic
in $T_{i}^{2}\times\{0\}$. Now since $\Sigma$ is incompressible, $\alpha$
must bound a disk $D$ in $\Sigma$, which is also a minimal surface embedded
in $(M^{4}(\tau_4),\gbar)$.

We claim that the minimal disk $D$ with $\partial{}D=\alpha$ must be
contained in $T_{i}^{2}\times[0,\tau_{4})$. 
Assume that $D$ is not entirely contained in $T_{i}^{2}\times(0,\tau_{4}])$.
Since the boundary of $T_{i}^{2}\times[0,\tau_{4}]$ is mean convex with
respect to both the inward normal vector and the modified metric $\gbar$,
according to the argument in \cite[pp. 155--156]{MY82b}, $\alpha$ bounds
an embedded least area minimal disk $D'\subset{}T_{i}^{2}\times[0,\tau_{4})$
(recall that $T_{i}^{2}\times\{\tau_{4}\}$ is {\tg} {\wrt} the metric $\gbar$)
and all of these kinds of least area minimal disks with the same boundary
$\alpha$ must be contained in $T_{i}^{2}\times[0,\tau_{4})$.
Recall that $D$ is assumed not to be entirely contained in 
$T_{i}^{2}\times(0,\tau_{4})$, so $\Area(D)>\Area(D')$, where the area
$\Area(\cdot)$ is with respect to the modified metric $\gbar$ on $M^{3}(\tau_4)$.
Since $\pi_{2}(M^{3}(\tau_4))=\pi_{2}(M^3)=0$, $D$ is isotopic to $D'$ with 
boundary fixed in $M^{3}(\tau_4)$. Let $\Pi$ be the surface defined by 
$\Pi=(\Sigma-D)\cup{}D'$, then $\Pi$ is isotopic to $\Sigma$ and 
$\Area(\Pi)<\Area(\Sigma)$, where the area $\Area(\cdot)$ is also with respect 
to the metric $\gbar$. 
But this contradicts the assumption that $\Sigma$ is a least area minimal
surface isotopic to $S$ in $(M^{3}(\tau_4),\gbar)$.
Therefore $D$ must be itself contained in $T_{i}^{2}\times[0,\tau_{4})$.

Next we choose $\varepsilon$ sufficiently small such that
$\alpha\times[0,\varepsilon]\subset\Sigma'$, therefore each simple
closed curve $\alpha_{t}=\Sigma'\cap{}(T^{2}_{i}\times\{t\})$ is null
homotopic in $T^{2}_{i}\times\{t\}$ for $0\leq{}t\leq\varepsilon$,
then similarly we have a minimal disk $D_{t}$ such that
$D_t\subset\Sigma\cap(T^{2}_{i}\times[t,\tau_{4}])$
for $0\leq{}t\leq\varepsilon$. But this is impossible since otherwise
$\Sigma$ would self-intersect (uncountably) infinitely many
times at $\alpha_t$ for $0<t\leq\varepsilon$.

Therefore, if $\Sigma'$ is is a surface of negative Euler characteristic in
$T_{i}^{2}\times[0,\tau_{4})$, then no closed essential curves in $\Sigma'$ is
null homotopic in $T_{i}^{2}\times[0,\tau_{4})$. But this is impossible by
applying Lemma~\ref{euler}. 

\vskip 0.1in
\textbf{Case II}: Suppose that $\Sigma'$ is an annulus such that at least
one of its boundary components, say $\alpha$, is null-homotopic in
$T_{i}^2\times\{0\}$. Then apply the
similar argument as in Case I (in this case, $\alpha$ is also a peripheral
essential curve in $\Sigma'$), we know that this is impossible.

Thus each component of
$\Sigma\cap\left(T^{2}_{i}\times[0,\tau_{4})\right)$ is either a
minimal disk whose boundary is a null-homotopic {\jc} in
$T^{2}_{i}\times\{0\}$, or a minimal annulus whose boundary consists
of two essential {\jc}s in $T^{2}_{i}\times\{0\}$.
It's easy to see that each component of
$\Sigma\cap(T^{2}_{i}\times[0,\tau_{4}))$ is boundary compressible.
\ep

\subsection{Good positioned {\jc}s on tori}
We start by making a definition of {\jc}s being in good position on a torus.
This will be important for what follows.

\begin{Def}\label{good}
Let $M^3$ be a {\chtm} and $\CC=T^{2}\times[0,\infty)$ be a maximal cusped region
of $M^3$. A {\jc} (i.e., simple closed) $\alpha\subset{}T^{2}\times\{\tau\}$ is
said to be in ``\emph{{\gp}}" if one of the lifts of $\alpha$ to $\H^3$
is contained in $A\times\{e^{\tau}\}$, where $A$ is the fundamental domain of
the parabolic group $\Gamma_{\infty}=\langle{}z\mapsto{}z+\mu,\ z\mapsto{}z+\nu\rangle$
in the horosphere $\{(x,y,1)\ |\ (x,y)\in\R^2\}$.
\end{Def}

From the above definition, we have the following statement:

\bpro\label{good2}
A {\jc} $\alpha\subset{}T^{2}\times\{\tau\}$ is {\igp} if the
Euclidean length of $\alpha$ is less
than $\min\{2|\mu|,2|\nu|, 2|\mu \pm \nu|\}$, while if it is not {\igp} 
then the Euclidean length of $\alpha$ is at least
$\min\{2|\mu|,2|\nu|, 2|\mu \pm \nu|\}$.

If $\alpha\subset{}T^{2}\times\{\tau\}$ is an
\emph{essential} {\jc}, then $\alpha$ is \emph{not} {\igp}.
\epro

Recall from \eqref{4c} that we have $4$ constants:
$\tau_{j}=\log(j\cdot{}L_{0})$ for $j=1,2,3,4$, where the constant
$L_{0}$ is defined in \eqref{maxh}. And these constants are ordered:
$\tau_4 > \tau_3 > \tau_2> \tau_1 > 0$. As in the previous subsection, we assume
$\Sigma\cap\left(T^{2}_{i}\times[0,\tau_{2}]\right)$ is non-empty. We first
observe the following fact:

\bpro\label{claim2}
Let $\Sigma'$ be a component of
$\Sigma\cap\left(T^{2}_{i}\times[0,\tau_{4})\right)$.
If there exists some $\tau\in[0,\tau_2]$, such that
$\Sigma'\cap\left(T^{2}_{i}\times\{\tau\}\right)$ consists
of {\jc}s {\igp}, then each component of $\Sigma'\cap\left(T^{2}_{i}\times\{\tau'\}\right)$
is also {\igp} for all $\tau'\in[\tau,\tau_2]$.
\epro

\bp
By Theorem \ref{kg}, we can lift $(M^{3}(\tau_{4}),\gbar)$ to
the truncated negatively curved space $(\Omega(\tau_{4}),\gbar)$
such that $T^{2}_{i}\times\{0\}$ is lifted to
the horizontal plane passing through the point $(0,0,1)$. Suppose
that the barycenter of the fundamental domain $A_{i}$ of the parabolic
group generated by $z\mapsto{}z+\mu_{i}$ and $z\mapsto{}z+\nu_{i}$
is the point $(0,0,1)$.

Suppose $D$ is a component of $\Sigma'\cap\left(T^{2}_{i}\times[\tau,\tau_{4})\right)$
such that $\partial{}D\subset{}T^{2}_{i}\times\{\tau\}$ is {\igp},
then by the arguments in Proposition ~\ref{claim1}, and
Proposition ~\ref{good2}, $D$ must be a disk and
$\partial{}D$ must be a null-homotopic {\jc} in $T^{2}_{i}\times\{\tau\}$. Let
$\widetilde{D}$ be a lift of $D$ such that
$\partial\widetilde{D}\subset{}A_{i}\times\{e^{\tau}\}$.

We define the following:
\be\label{B}
   \Bcal_i=A_{i}\times[e^{\tau},4L_{0}]\ .
\ene
We want to show that $\widetilde{D}$ must be contained in $\Bcal_i$.
In fact, it is a minimal disk such that
$\partial\widetilde{D}\subset{}A_{i}\times\{e^{\tau}\}$
is null-homotopic. Then we are left with very few cases:
\ben
\item The minimal disk $\widetilde{D}$ doesn't have any subdisk
      below the horizontal plane through the point $(0,0,e^{\tau})$,
      since $D\subset{}T^{2}_{i}\times[\tau,\tau_{4})$ by the assumption.
\item Since all vertical planes are {\tg} (see Theorem \ref{tg} and
      Remark ~\ref{vert}), the minimal disk $\widetilde{D}$ does not
      have any subdisk outside $\Bcal_i$ by Hopf's {\maxp}.
\een
Thus $\widetilde{D}$ must be contained in the domain $\Bcal_i$. This is
certainly true for the other lifts of $D$ which are given by
$\gamma(\widetilde{D})$ for $\gamma\in\Gamma$.
By definition, for $\tau'\geq\tau$, each component of
$\Sigma'\cap\left(T^{2}_{i}\times\{\tau'\}\right)$ is {\igp}.
\ep

As a corollary, and taking advantage of Theorem ~\ref{barrier}
that we can use hemispheres as barriers, we find:

\bcor\label{sigma'}
If there exists some $\tau\in[0,\tau_2]$ such that
$\Sigma'\cap\left(T^{2}_{i}\times\{\tau\}\right)$ consists of {\jc}s
{\igp}, then $\Sigma'$ is contained in $T^{2}_{i}\times[0,\tau_{3}]$,
i.e. $\Sigma'$ is a least area disk or annulus {\wrt} the {\hym}.
\ecor

\bp
Recall from \eqref{U} and \eqref{f2}, the modified metric $\gbar$ is
flat for $t > 4L_0$, and hyperbolic when $t < 3L_0$. For convenience,
we denote two new constants: $L_{3}=\sqrt{e^{2\tau}+(\frac{L_0}{2})^2}$ and
$L_{4}=\frac{\sqrt{65}}{2}L_{0}$. Since $\tau\leq\tau_{2}=\log(2L_{0})$,
so we have
\be\label{Ls}
    L_{3}\leq \frac{\sqrt{17}}{2} L_0<3L_{0}< 4L_0<L_4.
\ene
Therefore $A_i \times\{L_4\}$ is {\tg} {\wrt} the metric $\gbar$.

We consider the subregion $\Bcal'_i$ of $\Bcal_i$, which is defined by
\beq
   \Bcal'_i=\Bcal_i\cap\left\{\bigcup_{L_{3}\leq{}r
   \leq{}L_{4}}S_{+}^{2}(r)\right\}\ .
\eeq
By Theorem ~\ref{barrier}, the subregion $\Bcal'_i$ is foliated by the
non-concave spherical caps {\wrt} the downward {\nv}s. By the definition of
$L_{0}$ in \eqref{maxh}, the spherical cap $\Bcal_i\cap{}S_{+}^{2}(L_{3})$
lies above $A_{i}\times\{e^{\tau}\}$.

Recall from the proof of Proposition ~\ref{claim2} that $D$ is a
component of $\Sigma'\cap\left(T^{2}_{i}\times[\tau,\tau_{4}]\right)$
such that $\partial{}D\subset{}T^{2}_{i}\times\{\tau\}$ is {\igp}, and $\widetilde{D}$ is a lift of $D$ such that
$\partial\widetilde{D}\subset{}A_{i}\times\{e^{\tau}\}$. Therefore by the
{\maxp}, $\widetilde{D}$ is contained in $\Bcal_i$ and below the spherical
cap $\Bcal\cap{}S_{+}^{2}(L_{3})$. In other words, the Euclidean height
of $\widetilde{D}$ is at most $L_{3}$.

By \eqref{Ls}, we have
$\widetilde{D}\subset{}A_{i}\times[e^{\tau},3L_{0}]$. This is true
for other lifts of $D$ which are given by $\gamma(\widetilde{D})$,
for all $\gamma\in\Gamma$. Since the {\kg} preserves the metric
$\gbar$ (Theorem ~\ref{kg}), we have
$D\subset{}T^{2}_{i}\times[\tau,\tau_{3}]$, and therefore
\beq
   \Sigma'\subset{}\left(T^{2}_{i}\times[0,\tau]\right)\cup
                   \left(T^{2}_{i}\times[\tau,\tau_{3}]\right)
           =T^{2}_{i}\times[0,\tau_{3}]\ .
\eeq
The proof of the Corollary is complete.
\ep

\subsection{Completing the proof}
First we need a version of the {\caf} modified from that in
\cite[p.399]{CG06}. The proof of \eqref{formula1} in the following
Lemma \ref{coarea} can be found in \cite{Wan12}.

\bl\label{coarea}
If $M^3$ is a Riemannian {\tm} with nonempty boundary $\partial{}M^3$, and $F$ is a
component of $\partial{}M^3$ such that its $s$-neighborhood $\Nscr_{s}(F)\subset{}M^3$ is a
trivial {\nb} over itself. If $\Sigma_1\subset{}M^3$ is a surface such that
$\Sigma_1\cap\Nscr_{s}(F)\ne\emptyset$, then
\be\label{formula1}
   \Area(\Sigma_1\cap\Nscr_{s}(F))= \int_{0}^{s}\int_{\Sigma_1\cap\partial\Nscr_{\tau}(F)}
   \frac{1}{\cos\theta}\,dld\tau\ ,
\ene
where the angle $\theta$ is defined as follows: For any point
$q\in\Sigma_1$, set $\theta(q)$ to be the angle between the tangent
space to $\Sigma_1$ at $q$, and the radial geodesic which is
through $q$ {\rm(}emanating from $q${\rm)} and is perpendicular to $F$.
\el

To complete the proof of Theorem ~\ref{main}, we just need to find one
$\tau \in [0,\tau_2]$ satisfying the assumption in Proposition ~\ref{claim2}.
And we show this $\tau$ may be chosen as just $\tau_2$:

\bt\label{claim3}
Let $\Sigma'$ be a component of
$\Sigma\cap\left(T^{2}_{i}\times[0,\tau_{4})\right)$,
then any component of
$\Sigma'\cap\left(T^{2}_{i}\times\{\tau_{2}\}\right)$ is a
{\jc} {\igp}.
\et

\bp
Assume that $\Sigma'$ is a component of $\Sigma\cap\left(T^{2}_{i}\times[0,\tau_{4}]\right)$ such that
at least one component of
$\Sigma'\cap\left(T^{2}_{i}\times\{\tau_{2}\}\right)$ is not {\igp},
then by Proposition ~\ref{claim2}, for each $\tau\in[0,\tau_2]$, $\Sigma'\cap\left(T^{2}_{i}\times\{\tau\}\right)$
has at least one component that is not {\igp}.

By Proposition ~\ref{good2}, for all $\tau\in[0,\tau_2]$, we have:
\be
   \Length\left(\Sigma'\cap\left(T^{2}_{i}\times\{\tau\}\right)\right)\geq{}
   \min\{2|\mu_i|,2|\nu_i|, 2|\mu_i \pm \nu_i|\}e^{-\tau}\ .
\ene

To apply the {\caf} \eqref{formula1}, we choose $F=T^{2}_{i}\times\{0\}$,
and for each $\tau\in[0,\tau_{2}]$, we set
\be
   \Nscr_{\tau}(F)=\left\{p\in{}T^{2}_{i}\times[0,\tau_{2}]\ |\
   \dist(p,F)\leq\tau\right\}\ ,
\ene
where $\dist(\cdot,\cdot)$ is the hyperbolic distance function. Now we apply the
{\caf} \eqref{formula1} to find:
\begin{align*}
   \Area\left(\Sigma'\cap\left(T^{2}_{i}\times[\tau_{1},\tau_{2}]\right)\right)
     &=\int_{\tau_1}^{\tau_2}\int_{\Sigma'\cap\partial\Nscr_{\tau}(F)}
        \frac{1}{\cos\theta}\,dld\tau\\
     &\geq\int_{\tau_1}^{\tau_2}\Length(\Sigma'\cap\partial\Nscr_{\tau}(F))\,d\tau\\
     &\geq\int_{\tau_1}^{\tau_2}\min\left\{2|\mu_i|,2|\nu_i, 2|\mu_i \pm \nu_i||\right\}e^{-\tau}\,d\tau\\
     &= \frac{\min\left\{|\mu_i|,|\nu_i|, |\mu_i \pm \nu_i|\right\}}{L_0}\\
     &> \frac{\min\{|\mu_{i}|\cdot|\nu_{i}|, |\mu_i-\nu_i|\cdot|\mu_i+\nu_i|\}}{L_0^2}\\
     &=\min\{|\mu_{i}|\cdot|\nu_{i}|, |\mu_i-\nu_i|\cdot|\mu_i+\nu_i|\}e^{-2\tau_{1}}\\
      &\geq\Area\left(T^{2}_{i}\times\{\tau_{1}\}\right)\ .
\end{align*}
Here we used the fact that $L_0 \ge |\mu_i| + |\nu_i|$ (\eqref{maxh}) and
$\tau_j = \log(jL_0)$ for $j=1,2$.

By Proposition ~\ref{claim1}, $\Sigma'$ is either a {\la} disk or a
{\la} annulus (that is boundary compressible), so we may isotope
$\Sigma'\cap(T^{2}_{i}\times[\tau_{1},\tau_{4}])$ to a disk or an annulus
$A$ contained in $T^{2}_{i}\times\{\tau_{1}\}$ such that
$\partial{}A=\Sigma'\cap(T^{2}_{i}\times\{\tau_{1}\})$. Let $\Sigma''$ be a new
surface defined by
\begin{equation*}
   \Sigma''=\left(\Sigma'\cap(T^{2}_{i}\times[0,\tau_{1}))\right)
            \cup{}A\ .
\end{equation*}
Then $\partial\Sigma''=\partial\Sigma'$ and
$\Sigma''$ is isotopic to $\Sigma'$ with boundary fixed in
$\Sigma'\cap(T^{2}_{i}\times[0,\tau_{4}])$. By the above inequality,
we have $\Area(\Sigma'')<\Area(\Sigma')$, but
this contradicts the fact that $\Sigma'$ is a
least area minimal surface in the region $T_{i}^{2}\times[0,\tau_{4}]$.
Therefore any component of
$\Sigma'\cap\left(T^{2}_{i}\times\{\tau_{2}\}\right)$ is a {\jc} {\igp}, and then
any component of $\Sigma\cap\left(T^{2}_{i}\times\{\tau_{2}\}\right)$ is also {\igp}.
\ep

We may now complete the proof: 

\bp[\bf{Proof of Theorem~\ref{main}}]\label{page:main theorem}
We consider two cases:

\textbf{Case I}: $S$ is assumed to be embedded in $M^3$.
This will complete the proof for Corollary ~\ref{main2}.

According to both Theorem 5.1 and the remarks before Theorem 6.12
in \cite{HS88}, there is an embedded incompressible least area minimal
surface $\Sigma$ isotopic to $S$ in $(M^{3}(\tau_4),\gbar)$.
By Theorem \ref{claim3}, all components of
$\Sigma\cap\left(T^{2}_{i}\times\{\tau_{2}\}\right)$ are {\igp}, then by
Corollary ~\ref{sigma'}, each component of
$\Sigma\cap(T^{2}_{i}\times[0,\tau_{4}))$ is disjoint
from $T^{2}_{i}\times(\tau_{3},\tau_{4}]$. Therefore we have
\begin{equation*}
   \Sigma\cap(T^{2}_{i}\times[0,\tau_{4}))
   \subset{}T^{2}_{i}\times[0,\tau_{3}]\ ,
   \quad\text{for}\ i=1,\ldots,k\ ,
\end{equation*}
which implies that $\Sigma$ is a {\ms} {\wrt} the {\hym}.

Next we claim that the minimal surface $\Sigma\subset{}M^{3}(\tau_3)$
is a least area minimal surface isotopic to $S$ in the cusped hyperbolic
$3$-manifold $M^3$.

In fact, let $L_{0}'$ be an arbitrary real number
such that $L_{0}'\geq L_{0}$, and let $\tau_{j}'=\log(j\cdot L_{0}')$ for
$j=1,2,3,4$. Obviously $\tau_{j}'\geq\tau_j$, and so
$M^3(\tau_{j})\subset{}M^3(\tau_{j}')$ for $j=1,2,3,4$.
We can construct the truncated $3$-manifold
$M^3(\tau_{4}')$ with a modified metric $\gbar'$ as in
$\S$\ref{sec:truncated space}, i.e.
\begin{itemize}
  \item $\gbar'|M^{3}(\tau_{3}')$ is hyperbolic, and
  \item $\gbar'|(M^{3}(\tau_{4}')-M^{3}(\tau_{3}'))$
         is defined as in $\S$\ref{sec:truncated space}.
\end{itemize}
and similarly any least area minimal surface $\Sigma'$ isotopic to $S$ in
$(M^3(\tau_{4}'),\gbar')$ must be contained in $M^3(\tau_{3})$, so it
can be considered as a minimal surface isotopic to $S$ in
$(M^3(\tau_{4}),\gbar)$. But we know that $\Sigma$ is the least area
minimal surface isotopic to $S$ in $(M^3(\tau_{4}),\gbar)$, so we must have
$\Area(\Sigma')\geq\Area(\Sigma)$ with respect to the hyperbolic metric on
$M^3$ (since both $\Sigma$ and $\Sigma'$ are contained in $M^3(\tau_{3})$ and
$\gbar|M^3(\tau_{3})=\gbar'|M^3(\tau_{3})$ are both hyperbolic).
Let $L_{0}'\to\infty$, we know that $\Sigma$ is a least area minimal surface
isotopic to $S$ in $M^3$.

\textbf{Case II}: $S$ is only assumed to be immersed in $M^3$.

By Theorem \ref{thm:LERF}, we may lift $S$ to an embedded nonseparating
closed {\is} in a \emph{finite cover}
$\widetilde{M}^3$ of $M^3$. It's easy to see that $\widetilde{M}^3$
is also a cusped hyperbolic $3$-manifold. Suppose that $\widetilde{M}^3$
has $\ell$ maximal cusped regions
$\widetilde{T}^{2}_{1}\times[0,\infty),\ldots,
\widetilde{T}^{2}_{\ell}\times[0,\infty)$ such that each
parabolic group corresponding to the horosphere
$\widetilde{T}^{2}_{i}\times\{0\}$ is generated by
$\langle{}z\mapsto{}z+\widetilde{\mu}_i,\ z\mapsto{}z+\widetilde{\nu}_i\rangle$ for
$i=1,\ldots,\ell$. Note that $\ell\geq{}k$, where $k$ is the number of the maximal
cusped regions of $M^3$. We define
\be\label{maxh2}
   \widetilde{L}_{0}=\max\big\{e^{\tau_{0}},\,
   |\widetilde{\mu}_{1}|+|\widetilde{\nu}_{1}|,\,\ldots,\,
   |\widetilde{\mu}_{\ell}|+|\widetilde{\nu}_{\ell}|\big\}\ > 0,
\ene
where $\tau_{0}>0$ is the same number as in $\S$\ref{sec:maximal cusps}, i.e.
$\tau_{0}$ is the smallest number such that each {\mcr} is disjoint from any
other {\mcr}s of $M^3$ and $\widetilde{M}^3$ respectively. Obviously
$\widetilde{L}_{0}\geq{}L_{0}$, where $L_{0}$ is defined by \eqref{maxh}.
Similarly, we also define
\be\label{eq:depth}
   \widetilde{\tau}_j=\log(j\cdot\widetilde{L}_{0})
   \qquad j=1,2,3,4\ .
\ene
Then $\widetilde{\tau}_j\geq\tau_j$ for $j=1,2,3,4$.
We still need some notations. It's easy to verify that the submanifold
\begin{equation*}
   \widetilde{M}^{3}(0)=\widetilde{M}^{3}-
\bigcup_{i=1}^{\ell}(\widetilde{T}^{2}_{i}\times[0,\infty))
\end{equation*}
of $\widetilde{M}^{3}$ is a finite cover of the submanifold $M^{3}(0)$
of $M^3$ defined by \eqref{quotient} for $\tau=0$.
Just as we did in $\S$\ref{sec:truncated space}, we can define
the truncated manifold $\widetilde{M}^3(\widetilde{\tau}_4)$
and the modified metric $\widetilde{g}$ on
$\widetilde{M}^3(\widetilde{\tau}_4)$ such that
\begin{itemize}
  \item $\widetilde{g}|\widetilde{M}^{3}(\widetilde{\tau}_3)$ is hyperbolic, and
  \item $\widetilde{g}|(\widetilde{M}^{3}(\widetilde{\tau}_4)-
         \widetilde{M}^{3}(\widetilde{\tau}_3))$ is defined as in
         $\S$\ref{sec:truncated space}.
\end{itemize}
By the construction,
$\widetilde{M}^3(\widetilde{\tau}_4)\to{}M^3(\widetilde{\tau}_4)$
is a finite cover. Let $\gbar$ be the modified metric on
$M^3(\widetilde{\tau}_4)$ defined by the finite cover map such that
\begin{itemize}
  \item $\gbar|M^{3}(\widetilde{\tau}_3)$ is hyperbolic,
  \item $\gbar|(M^{3}(\widetilde{\tau}_4)-M^{3}(\widetilde{\tau}_3))$
        is defined as in $\S$\ref{sec:truncated space}, and
  \item the finite cover
        $(M^{3}(\widetilde{\tau}_4),\widetilde{g})\to
        (M^{3}(\widetilde{\tau}_4),\gbar)$ is a local isometry.
\end{itemize}
According to the definition the metric $\gbar$ on $M^{3}(\widetilde{\tau}_4)$,
particularly we know that $\gbar|M^{3}(\tau_3)$ is hyperbolic since
$\tau_3\leq\widetilde{\tau}_3$.

Since $S$ is assumed to be immersed in $M^3(\widetilde{\tau}_4)$, according to
both Theorem 5.3 and the remarks before Theorem 6.12 in \cite{HS88}, there exists
an immersed {\la} {\ms} $\Sigma$ homotopic to $S$ in
$M^{3}(\widetilde{\tau}_4)$ with respect to the modified metric $\gbar$.
It is also incompressible in $M^{3}(\widetilde{\tau}_4)$. Note that it
might not be minimal {\wrt} the {\hym}.
On the other hand, since $\pi_{1}(M^{3}(\widetilde{\tau}_4))=\pi_{1}(M^3)$
is LERF by Theorem \ref{thm:LERF}, we may lift $S$ to an embedded closed {\is} $\widetilde{S}$ in $\widetilde{M}^{3}(\widetilde{\tau}_4)$, and
lift $\Sigma$ to a (possibly only immersed) {\ms} $\widetilde\Sigma$
homotopic to $\widetilde{S}$ in $\widetilde{M}^{3}(\widetilde{\tau}_4)$
with respect to the modified metric $\widetilde{g}$.
Since $\Sigma$ is a least area minimal surface homotopic to $S$ in
$M^{3}(\widetilde{\tau}_4)$ with respect to the modified metric $\gbar$,
the minimal surface $\widetilde\Sigma$ is a least area minimal
surface homotopic to $\widetilde{S}$ in
$(\widetilde{M}^{3}(\widetilde{\tau}_4),\widetilde{g})$.
Then using both Theorem 5.1 and the remarks before Theorem 7.1 in
\cite{FHS83}, the minimal surface $\widetilde\Sigma$ is an
\emph{embedded} {\la} {\ms} isotopic to $\widetilde{S}$ in
$(\widetilde{M}^{3}(\widetilde{\tau}_4),\widetilde{g})$.

Now just as we did in Case I, we apply both Corollary ~\ref{sigma'}
and Theorem \ref{claim3} to $\widetilde\Sigma$ in the non-positively
curved manifold $(\widetilde{M}^{3}(\widetilde{\tau}_4),\widetilde{g})$,
then we have
\begin{equation*}
   \widetilde\Sigma\cap(\widetilde{T}^{2}_{i}\times[0,\widetilde{\tau}_{4}))
   \subset{}\widetilde{T}^{2}_{i}\times[0,\widetilde{\tau}_{3}]\ ,
   \quad\text{for}\ i=1,\ldots,\ell\ .
\end{equation*}
This means
$\widetilde\Sigma\subset(\widetilde{M}^{3}(\widetilde{\tau}_3),\widetilde{g})$,
so it is a least area minimal surface with respect to the hyperbolic metric.
Therefore we have $\Sigma\subset(M^{3}(\widetilde{\tau}_3),\gbar)$, i.e.,
$\Sigma$ is an immersed closed least area minimal surface homotopic to
$S$ in $M^{3}(\widetilde{\tau}_3)$ with respect to the hyperbolic metric,
and furthermore it is also a least area minimal surface homotopic to $S$
in the cusped hyperbolic $3$-manifold $M^3$ as we did in Case I.
\ep

\bibliographystyle{amsalpha}
\bibliography{cusp}
\end{document}